\documentclass[11pt,reqno]{article}
\usepackage{amssymb, amsmath, amsthm, amsfonts, amscd, epsfig}
\usepackage[english]{babel}
\usepackage{graphicx, epsfig}
\usepackage{color}

\newtheorem{theorem}{Theorem}[section]

\newtheorem{lemma}[theorem]{Lemma}

\newcommand{\ds}{\displaystyle}
\newcommand{\pf}{\noindent {\sl Proof}. \ }
\newcommand{\p}{\partial}

\newcommand{\eqnref}[1]{(\ref {#1})}

\newcommand{\Cbb}{\mathbb{C}}

\newcommand{\Rbb}{\mathbb{R}}

\newcommand{\Ecal}{\mathcal{E}}

\newcommand{\Lcal}{\mathcal{L}}

\def\Bx{{\bf x}}

\newcommand{\Ga}{\alpha}

\newcommand{\Gd}{\delta}
\newcommand{\Ge}{\epsilon}

\newcommand{\Gvf}{\varphi}

\newcommand{\Gk}{\kappa}

\newcommand{\Gl}{\lambda}

\newcommand{\Gm}{\mu}
\newcommand{\Gv}{\nu}

\newcommand{\Gs}{\sigma}

\newcommand{\GD}{\Delta}

\newcommand{\GG}{\Gamma}

\newcommand{\beq}{\begin{equation}}
\newcommand{\eeq}{\end{equation}}

\def\ol{\overline}
\newcommand{\hatna}{\widehat{\nabla}}

\numberwithin{equation}{section}
\numberwithin{figure}{section}

\begin{document}

\title{A proof of the Flaherty-Keller formula on the effective property of densely packed elastic composites\thanks{This work is supported by NRF 2016R1A2B4011304 and 2017R1A4A1014735}}

\author{Hyeonbae Kang\thanks{Department of Mathematics, Inha University, Incheon
    22212, S. Korea (hbkang@inha.ac.kr)} \and Sanghyeon Yu\thanks{Seminar for Applied Mathematics, ETH Z\"urich, R\"amistrasse 101, CH-8092 Z\"urich, Switzerland (sanghyeon.yu@sam.math.ethz.ch)}}

\maketitle

\begin{abstract}
We prove in a mathematically rigorous way the asymptotic formula of Flaherty and Keller on the effective property of densely packed periodic elastic composites with hard inclusions. The proof is based on the primal-dual variational principle, where the upper bound is derived by using the Keller-type test functions and the lower bound by singular functions made of nuclei of strain. Singular functions are solutions of the Lam\'e system and capture precisely singular behavior of the sress in the narrow region between two adjacent hard inclusions.
\end{abstract}




\section{Introduction}

In a two phase composite where inclusions are close to each other and the material property, such as conductivity and elastic moduli, of the inclusion is of high contrast with that of the matrix, the effective property of the composite becomes singular. Several asymptotic formula (as the distance between inclusions tends to $0$) for the effective properties have been found in diverse contexts. We review some of them which are related to present work.  For a more complete account of such results, we refer readers to section 10.10 of Milton's book \cite{Milton-book-01}.

In \cite{Keller-JAP-63}, Keller considered a square array of circular cylinders when cylinders are nearly touching to each other and have extreme conductivities (perfectly conducting or insulating). It is observed that the effective conductivity blows up in the nearly touching limit and this divergence is due to the fact that the electric current flux is concentrated in the narrow region between adjacent cylinders. In fact, he derived an asymptotic formula for the effective conductivity as the distance of two inclusions tends to $0$, and showed its validity numerically.

Berlyand and Kolpakov \cite{BK-ARMA-01} considered general non-periodic composites of circular inclusions and showed that its effective conductivity can be well approximated by a discrete network. In the course of such investigation, they proved Keller's asymptotic formula in a rigorous manner using the primal-dual variational principle. Remarkably their network approach has been extended to concentrated suspensions (Stokes system) and asymptotic formula for the effective property of periodic suspension has been derived \cite{BBP-SIMA-05, BGN-ARMA-09}.

Keller's formula for the effective conductivity has been extended to elastic composites by Flaherty and Keller \cite{FK-CPAM-73} (see \eqnref{FKformula} of this paper). They obtained asymptotic formulas for the effective elastic moduli of a rectangular array of cylinders in the nearly touching limit, when the cylinder is either a hard inclusion or a hole. However, a rigorous proof of their formulas is still missing to the best of our knowledge, and it is the purpose of this paper to provide a rigorous proof of the formula when the cylinder is a hard inclusion.

The proof of this paper is based on the primal-dual variational principle, which is a standard tool. The novelty of the proof lies in the construction of test functions, especially those for the dual principle. The test functions for the dual principle are constructed using singular functions introduced by authors in \cite{KY}. These functions are elaborated linear combinations of nuclei of strain and capture precisely stress concentration between two adjacent hard inclusions. Nuclei of strain are columns of the Kelvin matrix of the fundamental solution to the Lam\'e system and their variants. One important feature is that they are solutions of the Lam\'e system.

It is worth mentioning that in presence of two inclusions with extreme material property the gradient blows up in between two inclusions. Such gradient blow-up represents either stress concentration or field enhancement. Precise quantitative study on the gradient blow-up in various context has been important theme of active research in last ten years or so. We refer to \cite{KY} and references therein for such a development.

This paper is organized as follows. In the next section we set up the problem and introduce the Flaherty-Keller formula. In section \ref{sec:variation} the primal-dual variational principle is introduced with a short proof. In section \ref{sec:test} we introduce the Keller-type functions and singular functions as test functions for primal and dual variational principles, respectively. The last section is to prove the Flaherty-Keller formula.

\section{The Flaherty-Keller formula}\label{sec:FK}

\begin{figure}[ht!]
\begin{center}
\epsfig{figure=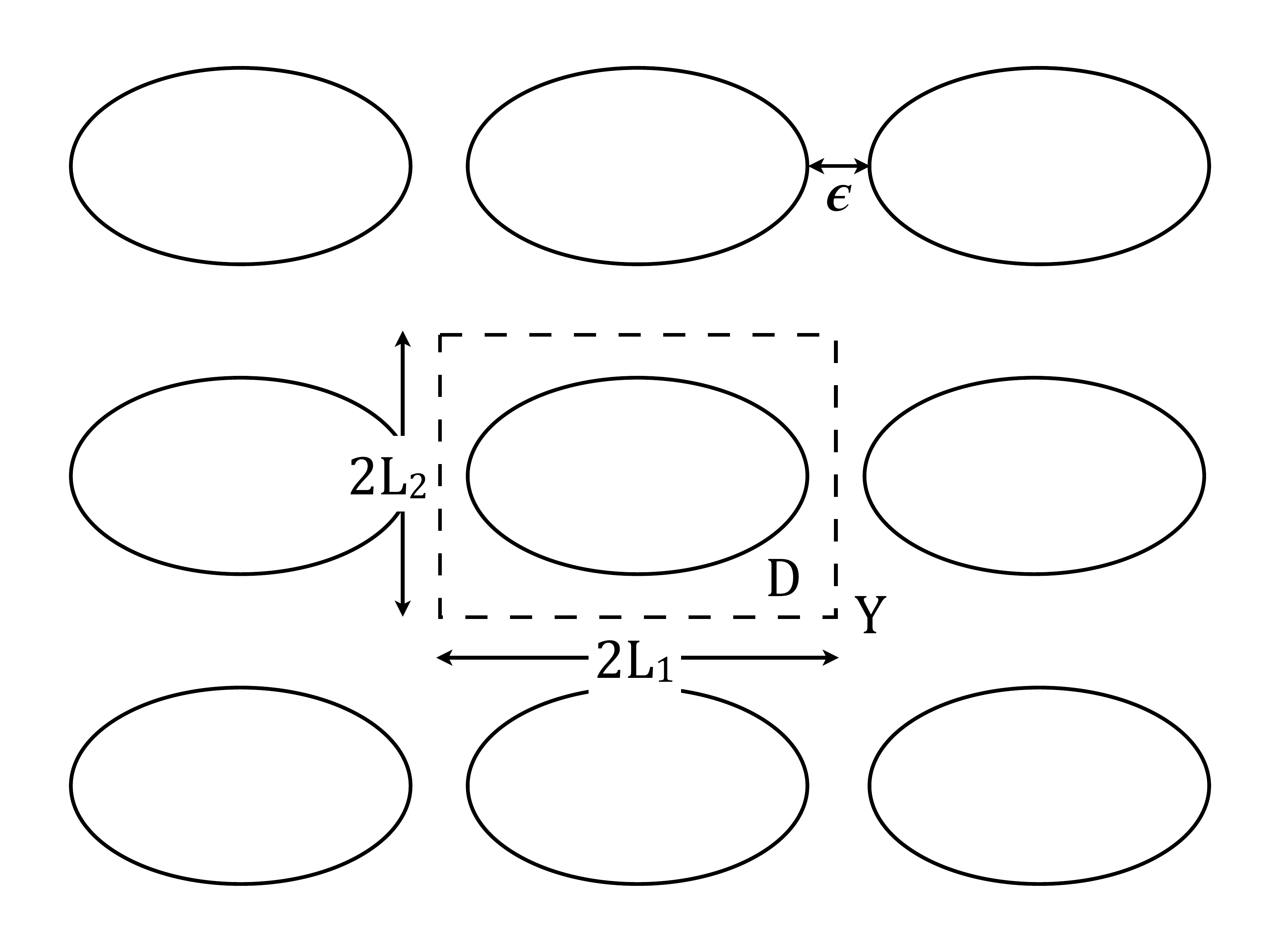,width=8cm}
\end{center}
\caption{Geometry of the composite}
\label{composite1}
\end{figure}

Let $L_1$ and $L_2$ be positive numbers and let $Y=(-L_1,L_1)\times(-L_2,L_2)$ denote the period cell.
Let $D \subset Y$ be a strictly convex domain containing the origin with the smooth boundary, which represents the two-dimensional cross section of the cylindrical inclusion.
We assume that $D$ is symmetric with respect to both $x$- and $y$-axes. Following \cite{FK-CPAM-73} we assume that the periodic inclusions are nearly touching in one direction and they are away from each other in the other direction. We assume that $D$ is close to the vertical boundary of $Y$, but away from the horizontal boundary. Let $\Ge/2$ be the distance between $D$ and the vertical boundary of $Y$, so that $\Ge$ becomes the distance between two adjacent inclusions. See Figure \ref{composite1}. It is worth mentioning that if $D$ is close to the horizontal boundary, then we can easily modify the Flaherty-Keller formula which is presented in Theorem \ref{thm:FK}.

Let $(\Gl,\Gm)$ be the pair of Lam\'{e} constants of $Y \setminus \overline{D}$ which satisfies the strong ellipticity conditions
$$
\Gm>0 \quad\mbox{and}\quad \Gl+\Gm>0.
$$
Then the elasticity tensor $\Cbb=(C_{ijkl})$ is given by
$$
C_{ijkl} = \Gl \Gd_{ij} \Gd_{kl} + \Gm (\Gd_{ik} \Gd_{jl} + \Gd_{il} \Gd_{jk}),
$$
where $\Gd_{ij}$ denotes Kronecker's delta. The Lam\'e operator $\Lcal_{\Gl,\Gm}$ of the linear isotropic elasticity is defined by
\beq
\Lcal_{\Gl,\Gm}  u:= \nabla \cdot \Cbb \hatna u = \Gm \GD  u+(\Gl+\Gm)\nabla \nabla \cdot  u,
\eeq
where $\hatna$ denotes the symmetric gradient, namely,
$$
\hatna  u = \frac{1}{2} \left( \nabla  u +\nabla  u^T\right) \quad (T \mbox{ for transpose}).
$$
The corresponding co-normal derivative $\p_\Gv u$ either on $\p D$ or on $\p Y$ is defined as
\beq
\p_\Gv u = (\Cbb \hatna u) n,
\eeq
where $n$ is the outward unit normal vector field to $\p D$ or $\p Y$, respectively.

Let
\beq\label{Psidef}
\Psi_1 =\begin{bmatrix} 1\\0 \end{bmatrix}, \quad
\Psi_2 =\begin{bmatrix} 0\\1 \end{bmatrix}.
\eeq
For $j=1,2$, let $v_j\in H^1(Y\setminus \overline{D})$ be the solution to the following problem:
\begin{align}
\begin{cases}
\ds
\Lcal_{\Gl,\Gm} v_j = 0, &\quad \mbox{ in }Y\setminus\overline{D}, \\
\ds v_j = 0, &\quad \mbox{ on }\p D,
\\[0.5em]
\ds v_j = \pm\frac{1}{2} \Psi_j, &\quad \mbox{ on } x=\pm L_1,
\\[0.5em]
\ds\p_\Gv v_j = 0, &\quad \mbox { on } y=\pm L_2.
\end{cases}
\end{align}
We extend $v_j$ to the whole space $\Rbb^2$ so that the extended function, denoted still by $v_j$, satisfies the following:
\beq
v_j (x+2L_1,y) = v_j(x,y)+\Psi_j, \quad v_j(x,y+2L_2) = v_j(x,y).
\eeq
In particular, $\hatna v_j$ is periodic. The effective extensional modulus $E_*$ and the effective shear modulus $\Gm_*$ are given by (see \cite{FK-CPAM-73})
\begin{align*}
E_* &= \frac{(1+\rho)(1-2\rho)}{1-\rho }\frac{L_1}{ L_2 }\int_{-L_2}^{L_2} \p_\Gv v_1 (L_1,y) \cdot \Psi_1 \,dy,
\\
\Gm_* &= \frac{L_1}{L_2}\int_{-L_2}^{L_2} \p_\Gv v_2 (L_1,y) \cdot \Psi_2 \,dy,
\end{align*}
where $\rho$ is Poisson's ratio, namely,
\beq\label{eqn_def_poisson}
\rho = \frac{\Gl}{2(\Gl+\Gm)}.
\eeq

Let $E$ be Young's modulus of the matrix, namely,
$$
E= \frac{\Gm(3\Gl+2\Gm)}{\Gl+\mu}.
$$
The asymptotic formulas for the effective elastic moduli obtained in \cite{FK-CPAM-73} are as follows:

\begin{theorem}[Flaherty-Keller formula]\label{thm:FK}
Let $\Gk_0$ be the curvature of $\p D$ at the closest point to the boundary $x=L_1$. The following hold:
\beq\label{FKformula}
E_* =  E\frac{L_1}{L_2}\frac{ \pi}{\sqrt{\Gk_0}}\frac{1}{\sqrt\Ge}+O(1)
\eeq
and
\beq\label{FKformula2}
\Gm_* = \Gm \frac{L_1}{ L_2 }\frac{ \pi}{\sqrt{\Gk_0}}\frac{1}{\sqrt\Ge}+O(1),
\eeq
as $\Ge \to 0$.
\end{theorem}

\section{Primal-dual variational principle}\label{sec:variation}

The effective moduli $E_*$ and $\Gm_*$ can be represented using energy integrals.
In fact, since $\hatna v_j$ is periodic and $n|_{x=-L_1}=-n|_{x=L_1}$, we see that
$$
\p_\Gv v_j (-L_1, y) = - \p_\Gv v_j (L_1, y).
$$
So, we have from the boundary condition of $v_j$ that
\begin{align*}
\int_{-L_2}^{L_2} \p_\Gv v_j(L_1,y) \cdot \Psi_j &= \int_{-L_2}^{L_2} \p_\Gv v_j(-L_1,y) \cdot (-\frac{1}{2} \Psi_j) + \int_{-L_2}^{L_2} \p_\Gv v_j(L_1,y) \cdot \frac{1}{2} \Psi_j \\
&= \int_{\p(Y\setminus \overline{D})} \p_\Gv v_j \cdot v_j = \int_{Y\setminus \ol{D}} \Cbb \hatna v_j: \hatna v_j.
\end{align*}
Here and throughout the paper, the expression $A:B$ for $2 \times 2$ matrices $A=(a_{ij})$ and $B=(b_{ij})$ indicates $\sum_{i,j} a_{ij} b_{ij}$.

Let
\beq
\Ecal_j := \int_{Y\setminus\overline{D}} \Cbb \hatna{v_j}:\hatna{v_j}.
\eeq
Then we have
\beq\label{eqn_eff_E_mu_energy}
E_* =  \frac{(1+\rho)(1-2\rho)L_1}{(1-\rho) L_2 }\Ecal_1
\quad\mbox{and}\quad
\Gm_* =  \frac{L_1}{L_2 }\Ecal_2.
\eeq

\begin{figure*}
\begin{center}
\epsfig{figure=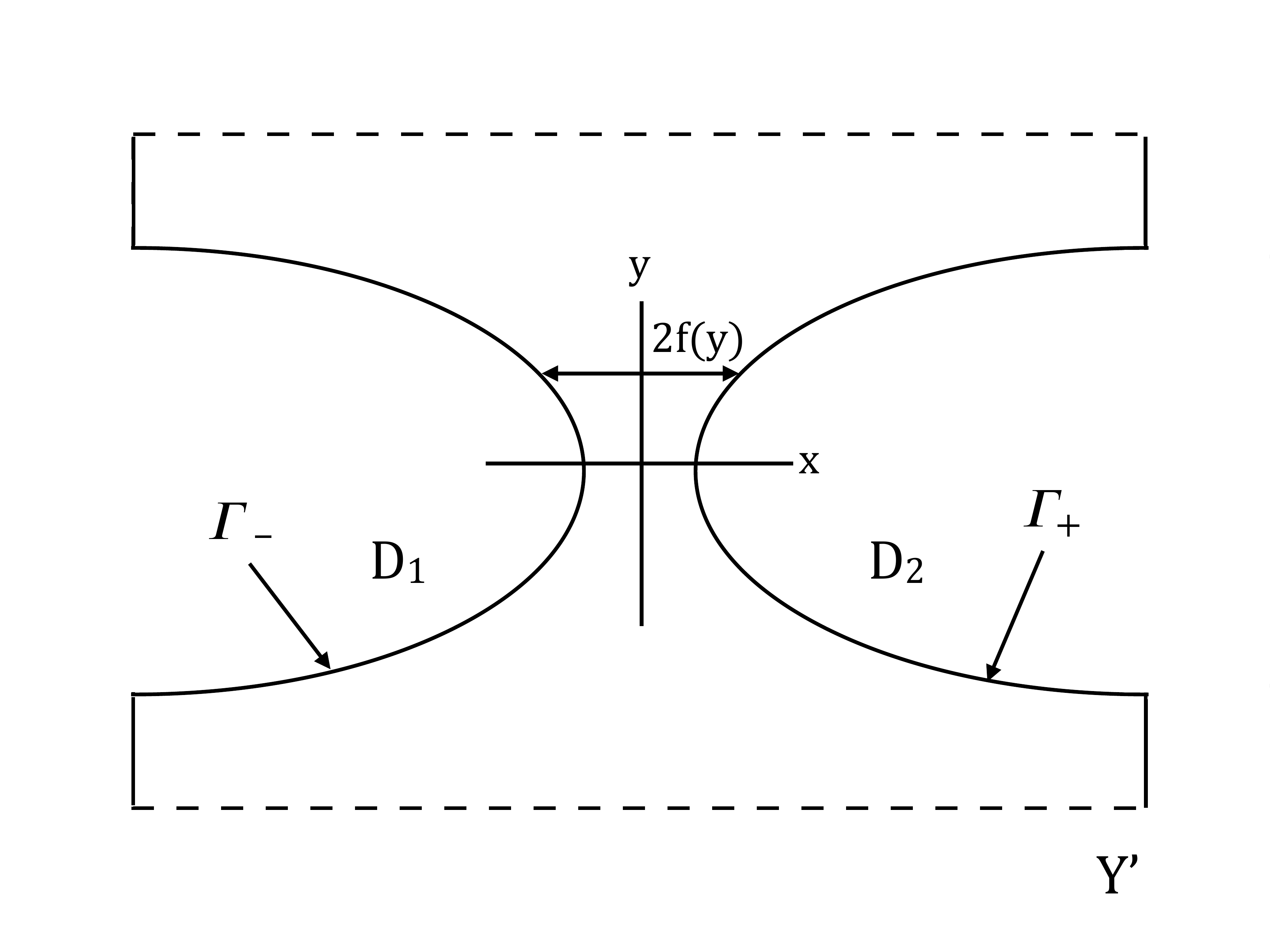,width=6cm}
\end{center}
\caption{Geometry of $Y'$}
\label{composite2}
\end{figure*}

It is more convenient to consider the energy integral $\Ecal_j$ in a translated cell $Y_t:=Y-(L_1,0)=(-2L_1, 0) \times (-L_2,L_2)$.
Let us denote $D_1 = D-(2L_1,0)$, $D_2=D$, and $Y'=Y_t \setminus\overline{D_1\cup D_2}$. Let
$\GG_-=(\p D_1 \cup \{x=-L_1\}) \cap \p Y'$ and $\GG_+=(\p D_2 \cup \{x=L_1\}) \cap \p Y'$. See Figure \ref{composite2} for the configuration of $Y'$.

Let us denote $v_j$ after translation by the same notation $v_j$. Then, by periodicity, we have
$$
\Ecal_j = \int_{Y'} \Cbb \hatna{v_j}:\hatna{v_j}.
$$
Note that, for $j=1,2$, the function $v_j|_{Y'} \in H^1(Y')$ is the solution to the following equation:
\beq
\begin{cases}
\ds
\Lcal_{\Gl,\Gm} v_j = 0, &\quad \mbox{ in }Y',
\\
\ds v_j = 0, &\quad \mbox{ on } \GG_-,
\\
\ds v_j = \Psi_j, &\quad \mbox{ on } \GG_+,
\\
\ds \p_\Gv v_j =0, &\quad \mbox{ on } \{y=\pm L_2\},
\end{cases}
\eeq

The following primal-dual variational principle is used:
\begin{lemma}[primal-dual variational principle]\label{primedual}
Let
\begin{align}
V_j&=\left\{\Gvf \in H^1(Y'):\Gvf|_{\GG_-}= 0, \,\Gvf|_{\GG_+} =\Psi_j\right\}, \label{Vj}
\\[0.2em]
W_j&=\left\{\Gs\in L^2_s(Y': \Rbb^{2\times 2}):\nabla\cdot\Gs=0, \,\Gs n=0 \mbox{ if } y=\pm L_2 \right\},
\end{align}
where $L^2_s(Y': \Rbb^{2\times 2})$ denotes the collection of the square integrable real symmetric $2 \times 2$ matrix valued functions. The following holds:
\begin{align}
\Ecal_{j} &= \min_{\Gvf\in V_j}\int_{Y'} \Cbb \hatna\Gvf:\hatna\Gvf
\label{eqn_primal_effective}
\\
&=\max_{\Gs\in W_j}
-\int_{Y'}  {\Gs}:\Cbb^{-1}{\Gs}+2\int_{\GG_+} \Gs n \cdot \Psi_j.
\label{eqn_dual_effective}
\end{align}
\end{lemma}

The primal principle \eqnref{eqn_primal_effective} is used to obtain the upper bound on $\Ecal_j$, while the dual principle \eqnref{eqn_dual_effective} is used for the lower bound. The primal-dual principle for the Laplace operator was used in \cite{BK-ARMA-01} to prove the Keller formula for the effective conductivity. The principle for the Lam\'e system may be well-known. However, we were not able to find a reference, and so we include a proof here.

\medskip

\noindent{\sl Proof of Lemma \ref{primedual}}.
By Green's identity and the boundary conditions in \eqnref{Vj}, we have
\begin{align*}
\int_{Y'}\Cbb \hatna v_j:\hatna \Gvf &= \int_{\p Y'} \p_\nu v_j\cdot \Gvf
=\int_{\GG_+}  \p_\nu v_j \cdot \Psi_j \quad\mbox{for all } \Gvf \in V_j.
\end{align*}
Therefore,
by the Cauchy-Schwartz inequality, we have
\begin{align*}
\int_{Y'} \Cbb\hatna v_j:\hatna v_j &=
\int_{Y'} \Cbb\hatna v_j:\hatna \Gvf
\\
&\leq \frac{1}{2} \left( \int_{D^e} \Cbb\hatna v_j:\hatna v_j
+ \int_{D^e} \Cbb\hatna \Gvf:\hatna \Gvf \right).
\end{align*}
This proves \eqnref{eqn_primal_effective}.

Note that $\Cbb \hatna v_j \in W_j$. If $\Gs \in W_j$,
by the divergence theorem and the fact that $\nabla \cdot \Gs=0$,
we have
\begin{align}
\int_{Y'} \Gs : \hatna v_j  &= -\int_{Y'} (\nabla\cdot \Gs) \cdot v_j + \int_{\p Y'}\Gs n \cdot  v_j
\notag
\\
&=\int_{|y|=L_2} 0 \cdot v_j
+ \int_{\GG_-} \Gs n \cdot 0
+ \int_{\GG_+} \Gs n \cdot \Psi_j =\int_{\GG_+} \Gs n \cdot \Psi_j.
\label{sigma_vj_int_by_parts}
\end{align}
By the Cauchy-Schwartz inequality, we have
$$
\int_{Y'} \Cbb \, f : g \leq \frac{1}{2} \Big(\int_{Y'} \Cbb \, f : f
+ \int_{Y'} \Cbb \, g: g \Big),
$$
for all $f,g\in L_s^2(Y':\mathbb{R}^{2\times 2})$.
So, by letting $f = \Cbb^{-1} \Gs$ and $g=\hatna v_j$, we obtain from \eqnref{sigma_vj_int_by_parts}
$$
\int_{\GG_+} \Gs n \cdot \Psi_j = \int_{Y'} \Gs : \hatna v_j \leq \frac{1}{2} \Big(\int_{Y'} \Gs : \Cbb^{-1} \Gs
+ \int_{Y'} \Cbb\hatna v_j : \hatna v_j\Big).
$$
Now, \eqnref{eqn_dual_effective} follows, and the proof is complete.
\qed

\section{Test functions}\label{sec:test}

Here we introduce test functions to be inserted into \eqnref{eqn_primal_effective} and \eqnref{eqn_dual_effective}. For that we first describe the geometry of two inclusions $D_1$ and $D_2$. Since $D_1$ and $D_2$ are strictly convex, there are unique points, one on $\p D_1$ and the other on $\p D_2$, which have the shortest distance. Denote them by $z_1\in \p D_1$ and $z_2\in \p D_2$. Let $\Gk_0$ be the common curvature of $\p D_j$ at $z_j$. Let $B_j$ be the disk osculating to $D_j$ at $z_j$ ($j=1,2$). Then the common radius $r_0$ of $B_j$ is given by  $r_0=1/\Gk_0$. Let $R_j$ be the reflection with respect to $\p B_j$ and let $p_1 \in B_1$ and $p_2 \in B_2$ be the unique fixed points of the combined reflections $R_1\circ R_2$ and $R_2\circ R_1$, respectively.
Then one can easily see that $p_1$ and $p_2$ are written as
\beq\label{pjdef}
p_1=(-a,0)\quad\mbox{and}\quad p_2=(a,0).
\eeq
where the constant $a$ is given by $a :=\sqrt{\Ge (4 r_0 + \Ge)}/2$,
from which one can infer
\beq\label{def_alpha}
a= \sqrt{\frac{\Ge}{\Gk_0}} + O(\Ge^{3/2}).
\eeq

Let us consider the narrow region between $D_1$ and $D_2$.
There exists $L>0$ independent of $\Ge$ and a strictly convex function $f:[-L,L]\rightarrow \Rbb$ such that
$z_2=(f(0),0)$, $f'(0)=0$, and $\p D_2$ is the graph of $f(y)$ for $|y|<L$, {\it i.e.}, $(f(y),y)\in \p D_2$.
Note that $\p D_1$ is the graph of $-f(y)$ for $|y|<L$ (see Figure \ref{composite2}).
Since $D_2$ is symmetric with respect to the $y$-axis, we have for $|y|<L$,
\beq\label{taylor}
f(y) = \frac{\Ge}{2} +\frac{\Gk_0}{2} y^2 + O(\Ge^2+y^4).
\eeq
We denote by $\Pi_L$ the narrow region between $D_1$ and $D_2$, namely,
\beq\label{narrowregion}
\Pi_L = \{ (x,y) \in\Rbb^2| -f(y)<x<f(y), \, |y|<L\}.
\eeq

For the primal problem \eqnref{eqn_primal_effective}, we use the following test function:
Let
\beq\label{keller}
\Gvf_j(x,y)=\frac{ x+f(y) }{2f(y)} \Psi_j, \quad (x,y)\in \Pi_{L}.
\eeq
Note that $\Gvf_j=0$ on $\GG_- \cap \p \Pi_L$ and $\Gvf_j=\Psi_j$ on $\GG_+ \cap \p \Pi_L$.
We then extend $\Gvf_j$ to $Y'$ so that $\Gvf_j|_{\GG_-}= 0$, $\Gvf_j|_{\GG_+}= \Psi_j$, $\|\Gvf_j \|_{H^1(Y'\setminus \Pi_L)} \le C$ for some $C$ independent of $\Ge$.
Then one can see that $\Gvf_j \in V_j$.

The function $(x+f(y))/2f(y)$ appearing in the definition of $\Gvf_j$ has been used in \cite{Keller-JAP-63} for derivation of the effective conductivity, and used in \cite{BK-ARMA-01} for its proof. For this reason, we call $\Gvf_j$ the Keller-type function. Recently the function $\Gvf_j$ was efficiently used by Bao {\it et al} \cite{BLL-ARMA-15, BLL-arXiv} to derive the upper bound on the blow-up rate of the gradient in presence of adjacent hard elastic inclusions. The upper bound in two dimensions turns out to be $\Ge^{-1/2}$ where $\Ge$ is the distance between two inclusions.

We emphasize that the function $\Gvf_j$ is not a solution of the Lam\'e system and does not seem to fit to the dual principle. In fact, we can modify the function $\Gvf_j$ so that the modified function becomes the solution of the Lam\'e system, and use it for the dual principle. But it does not yield the correct lower bound. In this paper we use singular functions introduced in \cite{KY} as test functions for the dual principle. Let $\GG = \left( \GG_{ij} \right)_{i, j = 1}^2$ be the Kelvin matrix of fundamental solutions to the Lam\'{e} operator $\Lcal_{\Gl, \Gm}$, namely,
\beq\label{Kelvin}
  \GG_{ij}(x) = \ds \Ga_1 \Gd_{ij} \ln{|\Bx|} - \Ga_2 \frac{x_i x_j}{|x| ^2},
\eeq
where
\beq
  \Ga_1 = \frac{1}{4\pi} \left( \frac{1}{\Gm} + \frac{1}{2 \Gm + \Gl} \right) \quad\mbox{and}\quad
  \Ga_2 = \frac{1}{4\pi} \left( \frac{1}{\Gm} - \frac{1}{2 \Gm + \Gl} \right).
\eeq
Singular functions are defined using the following functions as basic building blocks:
\beq\label{nuclei}
\GG(x)e_1, \quad \GG(x)e_2, \quad
\frac{x}{|x|^2}, \quad \frac{x^\perp}{|x|^2},
\eeq
where $\{e_1,e_2\}$ is the standard orthonormal basis in Cartesian coordinates and $x^\perp=(-x_2,x_1)$ if $x=(x_1,x_2)$.
These functions are known as nuclei of strain \cite{Love}. Singular functions are defined as follows:
\beq\label{Bqone} 
q_1(x) := \GG(x-p_1)e_1-\GG(x-p_2) e_1 + {\Ga_2 a} \left( \frac{x-p_1}{|x-p_1|^2}+\frac{x-p_2}{|x-p_2|^2}\right), 
\eeq
and 
\beq\label{Bqtwo}
q_2(x) :=  \GG(x-p_1)e_2 -\GG(x-p_2)  e_2 - {\Ga_2 a} \left( \frac{(x-p_1)^{\perp}}{|x-p_1|^2}+\frac{(x-p_2)^{\perp}}{|x-p_2|^2}\right),
\eeq
where $a$ is the number appearing in \eqnref{pjdef}.

It is shown in \cite{KY} that functions $q_1$ and $q_2$ capture the singular behavior of the gradient in presence of adjacent hard elastic inclusions. As a consequence it is proved that the upper bound $\Ge^{-1/2}$ mentioned above is actually the optimal bound on the gradient blow-up. We emphasize that $q_1$ and $q_2$ are solutions to the Lam\'e system.

To construct test functions in $W_j$, let
\beq\label{m2def}
m_1 := \frac{\pi(\Gl+2\Gm)}{\sqrt{\Gk_0}} \quad\mbox{and}\quad
m_2 := \frac{\pi \Gm}{\sqrt{\Gk_0}},
\eeq
and let
\beq\label{psijS}
\Gs_j^S :=\frac{m_j}{\sqrt\Ge}\Cbb\hatna q_j.
\eeq
Since $q_j$ is a solution of the Lam\'{e} system, $\Gs_j^S$ satisfies $\nabla \cdot \Gs_j^S =0$. But the function $\Gs_j^S$ does not belong to $W_j$ because
\beq\label{notzero}
\Gs_j^S n|_{y=\pm L_2}=\pm \Gs_j^S e_2\neq 0.
\eeq
To see this, recall that, for a displacement field $u$, its associated stress tensor $\Gs=\mathbb{C}\widehat\nabla u$ is represented by
$\Gs=(\Gs_{ij})$ where
\begin{align*}
\Gs_{11} &= (\Gl+2\mu) \p_1 u_1 + \Gl \p_2 u_2,
\\
\Gs_{22} &= \Gl \p_1 u_1 +  (\Gl+2\mu) \p_2 u_2,
\\
\Gs_{12} &= \Gs_{21}=\mu (\p_2 u_1 + \p_1 u_2).
\end{align*}
So we obtain
\beq\label{bounded}
\Gs_j^S e_2 = \frac{m_j}{\sqrt\Ge} \big( \mu (\p_2 q_{j1} + \p_1 q_{j2} ),\ \Gl \p_1 q_{j1} +  (\Gl+2\mu) \p_2 q_{j2} \big)^T,
\eeq
from which one can easily see that \eqnref{notzero} holds.

We now modify $\Gs_j^S$ by adding a function. Let
\begin{align}
F_j(x,y)&= -\frac{y+L_2}{2L_2} \big[ \Gs_j^S(x,L_2) + \Gs_j^S(x,-L_2) \big] e_2 + \Gs_j^S(x,-L_2)e_2,  \label{Fj}
\\
G_j(x) &= \frac{1}{2L_2} \int^x_0
\big[ \Gs_j^S(x',L_2) + \Gs_j^S(x',-L_2)\big] e_2
 \,dx'.
\end{align}
Let $\Gs_j^c$ be the $2\times 2$ matrix-valued function having $G_j$ and $F_j$ as its columns, namely,
$$
\Gs_j^c(x,y)=\begin{bmatrix}
G_j(x) & F_j(x,y)
\end{bmatrix}, \quad (x,y) \in Y'.
$$
Then, one can check that
\beq\label{eqn_psic_properties}
\nabla \cdot \Gs_j^c =0, \quad \Gs_j^c n|_{y=\pm L_2} = -\Gs_j^S n|_{y=\pm L_2}.
\eeq
In fact, since
\begin{align*}
\p_1(\Gs_j^c e_1) &= \p_1 G_j(x,y) =\frac{1}{2L_2} \big[ \Gs_j^S(x,L_2) + \Gs_j^S(x,-L_2)\big] e_2,
\\
\p_2(\Gs_j^c e_2) &= \p_2 F_j(x,y)= -\frac{1}{2L_2} \big[ \Gs_j^S(x,L_2) + \Gs_j^S(x,-L_2)\big] e_2,
\end{align*}
we have
$$
\nabla \cdot \Gs_j^c = \p_1(\Gs_j^c e_1) + \p_2(\Gs_j^c e_2) = 0.
$$
If $y= \pm L_2$, we see from \eqnref{Fj} that
$$
\Gs_j^c n = \Gs_j^c e_2  = F_j (x,\pm L_2) =\mp \Gs_j^S(x,\pm L_2) e_2 = - \Gs_j^S n .
$$

Let
\beq\label{dualfunc}
\Gs_j:=\Gs_j^S + \Gs_j^c.
\eeq
Then, from \eqnref{eqn_psic_properties} and the fact that $\nabla \cdot \Gs_j^S =0$, we easily see that $\Gs_j \in W_j$.

The following evaluations of integrals are obtained in \cite{KY} (Lemma 4.7 and 4.8):
\begin{align}
\int_{\p D_i} \p_\Gv q_j \cdot \Psi_k &=(-1)^{i}\Gd_{jk} , \quad i, j, k=1,2, \label{qjpsik} \\
\int_{\p D_1 \cup \p D_2} \p_\Gv q_j \cdot q_j &= {m_j^{-1}}\sqrt\Ge + O(\Ge), \quad j=1,2.
\label{qjqj}
\end{align}
We emphasize the signs in the above are opposite to those in Lemma 4.7 and 4.8 in \cite{KY} since the normal vector $n$ in $\p_\Gv q_j= (\Cbb \hatna q_j) n$ in this paper directed outward to $\p Y'$ (and hence inward to $\p D_i$). We also invoke an estimate from \cite{KY} (Lemma 3.4): 
\beq\label{lem34}
|q_j(x)| + |\nabla q_{j}(x)| \lesssim \sqrt\Ge \quad\mbox{for all } x \in \Rbb^2 \setminus (D_1 \cup D_2 \cup \Pi_{L}) , \quad j=1,2.
\eeq
Here and afterwards, the expression $A \lesssim B$ implies that there is a constant $C$ independent of $\Ge$ such that $A \le CB$. 

One can see from \eqnref{def_alpha}, \eqnref{bounded} and \eqnref{lem34} that $|\Gs_j^S n|_{y=\pm L_2}|\lesssim 1$. So we have
\beq\label{eqn_Bpsic_bdd}
|\Gs_j^c(x,y)| \lesssim 1, \quad \mbox{for } (x,y) \in Y'.
\eeq

\section{Proof of Theorem \ref{thm:FK}}\label{sec:proof}

We first derive an upper bound using the primal principle \eqnref{eqn_primal_effective}.
\begin{lemma}[Upper bound]\label{lem_eff_upper}
We have
\beq\label{upper}
\Ecal_1 \leq (\Gl+2\Gm)\frac{ \pi}{\sqrt{\Gk_0}}\frac{1}{\sqrt\Ge}+O(1)
\eeq
and
\beq\label{upper2}
\Ecal_2 \leq \Gm\frac{ \pi}{\sqrt{\Gk_0}}\frac{1}{\sqrt\Ge}+O(1).
\eeq
\end{lemma}
\pf
To prove \eqnref{upper}, we use \eqnref{eqn_primal_effective} with $\Gvf=\Gvf_1$, where $\Gvf_1$ is the function defined with \eqnref{keller}. Note that since $\|\Gvf_j \|_{H^1(Y'\setminus \Pi_L)} \le C$ for some $C$ independent of $\Ge$, it suffices to estimate $\int_{\Pi_L}\Cbb \hatna{\Gvf_1}:\hatna{\Gvf_1}$.

Let $\p_1$ and $\p_2$ denote partial derivatives with respect to $x$ and $y$ variables, respectively. Straightforward computations yield
\begin{align}
\Cbb \hatna{\Gvf_1}:\hatna{\Gvf_1}
&=
\begin{bmatrix}
(\Gl+2\Gm)\frac{1}{2f(y)} & \Gm \p_2 \frac{x}{2f(y)}
\\[0.5em]
\Gm \p_2 \frac{x}{2f(y)} & 0
\end{bmatrix}
:
\begin{bmatrix}
\frac{1}{2f(y)} & \p_2 \frac{x}{4f(y)}
\\[0.5em]
\p_2 \frac{x}{4f(y)} & 0
\end{bmatrix}
\nonumber
\\
&=\frac{\Gl+2\Gm}{4} \frac{1}{f(y)^2}+ \frac{\Gm}{4} \frac{x^2 f'(y)^2}{f(y)^4}.
\end{align}
So we have
\begin{align}
\int_{\Pi_L}\Cbb \hatna{\Gvf_1}:\hatna{\Gvf_1}  &=\frac{\Gl+2\Gm}{4} \int_{\Pi_L} \frac{1}{f(y)^2} + \frac{\Gm}{2} \int_{\Pi_L}
\frac{x^2 f'(y)^2}{f(y)^4} =: I+II.
\end{align}
Thanks to the Taylor expansion \eqnref{taylor} of $f$, we have
\begin{align*}
\frac{4}{\Gl+2\Gm}I&= \int_{-L}^L \int_{-f(y)}^{f(y)} \frac{1}{f(y)^2} dxdy = \int_{-L}^L \frac{2}{f(y)} dy \\
&= \int_{-L}^L \frac{2}{(\Ge+\Gk_0 y^2)/2 } dy + \int_{-L}^L \frac{2}{f(y)}-\frac{2}{(\Ge+\Gk_0 y^2)/2 } dy
\\
&= \int_{-\infty}^\infty \frac{2}{(\Ge+\Gk_0 y^2)/2 } dy  +O(1)+\int_{-L}^L \frac{O(y^4)}{f(y) (\Ge+\Gk_0 y^2)/2} dy
\\
&= \frac{4\pi}{\sqrt\Gk_0 \sqrt{\Ge}} +O(1).
\end{align*}
We also have 
\begin{align}
|II|&\lesssim \int_{\Pi_L} \frac{(\Ge+y^2)^2 y^2}{(\Ge+y^2)^4}
\lesssim \int_{-L}^L \frac{ y^2}{\Ge+y^2} dy \lesssim 1.
\end{align}
Therefore we obtain
$$
\int_{\Pi_L}\Cbb \hatna{\Gvf_1}:\hatna{\Gvf_1}
= \frac{ \pi(\Gl+2\Gm)}{\sqrt{\Gk_0}}\frac{1}{\sqrt\Ge} +O(1).
$$
So, \eqnref{upper} follows from the primal principle \eqnref{eqn_primal_effective}.

One can prove \eqnref{upper2} in the exactly same manner using $\Gvf_2$ defined in \eqnref{keller}.
\qed

\medskip 

We now derive a lower bound using the dual variational princicple \eqnref{eqn_dual_effective}.
\begin{lemma}[Lower bound]\label{lem_eff_lower}
We have
\beq\label{lower}
\Ecal_1 \geq (\Gl+2\Gm)\frac{ \pi}{\sqrt{\Gk_0}}\frac{1}{\sqrt\Ge}+O(1)
\eeq
and
\beq\label{lower2}
\Ecal_2 \geq \Gm\frac{ \pi}{\sqrt{\Gk_0}}\frac{1}{\sqrt\Ge}+O(1).
\eeq
\end{lemma}
\pf
Let $\Gs_j$ be the function in $W_j$ defied by \eqnref{dualfunc}. Then we have
\begin{align*}
& -\int_{Y'}  {\Gs_j}: \Cbb^{-1}{\Gs_j}+2\int_{\GG_+} \Gs_j n \cdot \Psi_j\\
&= \bigg[-\int_{Y'}  {\Gs_j^S}:\Cbb^{-1}{\Gs_j^S}+2\int_{\GG_+} \Gs_j^S n \cdot \Psi_j \bigg]
+\bigg[-\int_{Y'}  {\Gs_j^c}:\Cbb^{-1}{\Gs_j^c}+2\int_{\GG_+} \Gs_j^c n \cdot \Psi_j\bigg] \nonumber
\\
&=: I_j +II_j.
\end{align*}
From \eqnref{eqn_Bpsic_bdd}, it is clear that
$
|II_j| \lesssim 1
$.

Now we estimate $I_j$. From the definition \eqnref{psijS} of $\Gs_j^S$, we have
\beq\label{Ij}
I_j =-\frac{m_j^2}{\Ge}\int_{Y'} \Cbb \hatna q_j:\hatna q_j
+ \frac{2m_j}{\sqrt\Ge}\int_{\GG_+} \p_\Gv q_j \cdot \Psi_j .
\eeq
Since $q_j$ is a solution of the Lam\'e system, we obtain by the divergence theorem that
$$
\int_{Y'} \Cbb \hatna q_j:\hatna q_j = \int_{\p Y'} \p_\Gv q_j \cdot q_j.
$$
Let $A:= (\p D_1 \cup \p D_2) \setminus \p Y'$ and $B:= \p Y' \setminus (\p D_1 \cup \p D_2)$ so that $\p Y' = (\p D_1 \cup \p D_2) \setminus A \cup B$, and
$$
\int_{\p Y'} \p_\Gv q_j \cdot q_j = \int_{\p D_1 \cup \p D_2} - \int_A + \int_B \p_\Gv q_j \cdot q_j.
$$
Since $A$ and $B$ are away from $\Pi_L$, we infer from \eqnref{lem34} that
$$
\left| - \int_A + \int_B \p_\Gv q_j \cdot q_j \right| \lesssim \Ge.
$$
It then follows from \eqnref{qjqj} that
\beq\label{intone}
\int_{Y'} \Cbb \hatna q_j:\hatna q_j = m_j^{-1} \sqrt{\Ge} + O(\Ge).
\eeq

Similarly, we write $\GG_+ = (\p D_2 \cup (\GG_+ \setminus \p D_2)) \setminus (\p D_2 \setminus \GG_+)$. Since $\GG_+ \setminus \p D_2$ and $\p D_2 \setminus \GG_+$ are away from $\Pi_L$, we infer that
$$
\int_{\GG_+} \p_\Gv q_j \cdot \Psi_j = \int_{\p D_2} \p_\Gv q_j \cdot \Psi_j + O(\sqrt{\Ge}) .
$$
So, \eqnref{qjpsik} yields
\beq\label{inttwo}
\int_{\GG_+} \p_\Gv q_j \cdot \Psi_j = 1+ O(\sqrt{\Ge}).
\eeq

It then follows from \eqnref{Ij}, \eqnref{intone} and \eqnref{inttwo} that
$$
I_j = m_j \sqrt{\Ge} + O(1).
$$
Now \eqnref{lower} and \eqnref{lower2} follow from \eqnref{m2def}. This completes the proof.
\qed

\medskip

Lemma \ref{lem_eff_upper} and \ref{lem_eff_lower} certainly lead us to Theorem \ref{thm:FK}.


\end{document}